\documentclass[3p,12pt]{elsarticle}

\usepackage{amssymb,amsmath,amsthm}
\newtheorem{them}{Theorem}[section]
\newtheorem{lem}{Lemma}[section]
\newtheorem{coro}{Corollary}[section]
\newtheorem{pro}{Proposition}[section]
\newtheorem{definition}{Definition}[section]

\begin{document}

\begin{frontmatter}




\title{On Riemann Integration in Metrizable Vector Spaces\tnoteref{1}}
\tnotetext[1]{This research was supported by the National Natural
Science Foundation of China (Grant Nos. 11826204, 11826206 and 11771384), the Natural Science Foundation of Yunnan Province of China (No. 2018FB004) and the Scientific Research Foundation of Yunnan University under grant No. 2018YDJQ010, and by Joint Key Project of Yunnan Provincial Science and Technology Department and Yunnan University (No. 2018FY001(-014)) and IRTSTYN.}

\author[2]{Zhou Wei\corref{4}}
\ead{wzhou@ynu.edu.cn}
\author[2]{Zhichun Yang}
\ead{yzhichun@mail.ynu.edu.cn}
\author[3]{Jen-Chih Yao}
\ead{yaojc@math.nsysu.edu.tw}
\cortext[4]{Corresponding author.}
\address[2]{Department of Mathematics, Yunnan University, Kunming 650091, People's Republic of China}

\address[3]{Department of Mathematics, Zhejiang Normal University, Jinhua, 321004, P. R. China  \\ and Research Center for Interneural Computing, China Medical University Hospital,  China Medical University, Taichung 40402, Taiwan}

\begin{abstract}
In classical analysis, Lebesgue first proved that $\mathbb{R}$ has the property that each Riemann integrable function from $[a,b]$ into $\mathbb{R}$ is continuous almost everywhere. This property is named as the Lebesgue property. Though the Lebesgue property may be breakdown in many infinite dimensional spaces including Banach or quasi Banach spaces, to determine spaces having this property is still an interesting problem. In this paper, we study Riemann integration for vector-value functions in metrizable vector spaces and prove the fundamental theorems of calculus and primitives for continuous functions. Further we discovery that $\mathbb{R}^{\omega}$, the countable infinite product of $\mathbb{R}$ with itself equipped with the product topology, is a metrizable vector space having the Lebesgue property  and prove that $l^p(1<p\leq+\infty)$, as subspaces of $\mathbb{R}^{\omega}$, possess the Lebesgue property although they are Banach spaces having no such property.
\end{abstract}


\end{frontmatter}

\section{Introduction}
The study on Riemann integration of functions mapping a closed interval into a Banach space was first studied by Graves \cite{Gr}. We refer the readers to bibliographies \cite{A,Go,Ha,Ro} for more details on Riemann integration. One interesting problem concerning the Riemann integration is to determine spaces having the Lebesgue property; i.e., every Riemann integrable function is continuous almost everywhere. In 1970s, Nemirovski, Ochan and Rejouani \cite{NOR} and da Rocha \cite{Ro} showed that $l^1$ has the Lebesgue property. In 1980s Haydon \cite{Ha} proved that for a stable Banach space with uniformly separable types, the Lebesgue property is equivalent to the Schur property (that is each weakly null sequence converges in norm). In 1990s, Gordon \cite{Go} listed many classical Banach spaces including $l^p(1<p<+\infty)$, $c_0$ and $L^p(1\leq p<+\infty)$ that do not have the Lebesgue property and published a non-classical space having the Lebesgue property; the Tsirelson space $T$. In 2008, Naralenkov \cite{Na} extended Gordon's results to prove that an asymptotic $l^1$ Banach space has the Lebesgue property.

The main work of this paper is to study Riemann integration in the context of metrizable vector spaces. Many of the real-value results concerning Riemann integration remain valid in this vector case. In particular we provide the fundamental theorems of calculus and primitives for continuous functions. As main results of this paper, we prove that two  metrizable vector spaces of $l^1(\Gamma)$ ($\Gamma$ uncountable) and $\mathbb{R}^{\omega}$ equipped with the product topology have the Lebesgue property. Further we prove that $l^p(p\geq1)$, as subspaces of $\mathbb{R}^{\omega}$, possess the Lebesgue property though $l^p(1<p\leq+\infty)$ are Banach spaces having no such property.

The remainder of this paper is organized as follows. Section 2 is devoted to preliminaries used in this paper. We consider the Riemann integration of functions in metrizable vector spaces. Several properties on the Riemann integration are provided. Two fundamental theorems of calculus and primitives on Riemann integration are proved. Section 3 contains main work of this paper. We show two metrizable vector spaces having the Lebesgue property; that are $l^1(\Gamma)$ ($\Gamma$ uncountable) and $\mathbb{R}^{\omega}$ equipped with product topology. The conclusion of this paper is given in Section 4.

\section{Preliminaries}
Let $X$ be a vector space over the real field $\mathbb{R}$ and it is equipped with a metric $\bf d$. The metric $\bf d$ is complete if every Cauchy sequence in $X$ converges to some point in $X$, and $\bf d$ is invariant if ${\bf d}(x+z, y+z)={\bf d}(x,y)$ for any $x,y,z\in X$. A subset set $\Omega$ in $ X$ is bounded if there exists $L\in (0,+\infty)$ such that ${\bf d}(0,x)\leq L$ for all $x\in \Omega$.

\begin{definition}
Let $[a,b]$ be an closed interval. A partition of $[a,b]$ is a finite set of points $\{t_i:0\leq i\leq N\}$ that satisfy $a=t_0<t_1<\cdots<t_{N-1}<t_N=b$. A tagged partition of $[a,b]$ is a partition $\{t_i:0\leq i\leq N\}$ of $[a,b]$ together with a set of points $\{s_i:1\leq i\leq N\}$ satisfying $s_i\in[t_{i-1}, t_i]$ for each $i$.

Let $\Delta=\{(s_i, [t_{i-1}, t_i]):1\leq i\leq N\}$ be a tagged partition of $[a,b]$. The points $\{t_i:0\leq i\leq N\}$ are said to be the points of the partition, the intervals $\{ [t_{i-1}, t_i]:1\leq i\leq N\}$ are said to be the intervals of the partition and the points $\{t_i:0\leq i\leq N\}$ are said to be the tags of the partition. Denote $|\Delta|:=\max_{1\leq i\leq N}(t_i-t_{i-1})$.

Let $\Delta_1$ and $\Delta_2$ be tagged partitions of $[a,b]$. The tagged partition $\Delta_1$ is a refinement of the tagged partition $\Delta_2$ if the points of $\Delta_2$ form a subset of the points of $\Delta_1$. In this case, we say that $\Delta_1$ refines $\Delta_2$.
\end{definition}

\begin{definition}
Let $X$ be a vector space equipped with a complete and invariant metric $\bf d$ and $f:[a,b]\rightarrow X$. The function $f$ is said to be Riemann integrable on $[a,b]$ if there exists a vector $\bar x\in X$ with the following property: for any $\varepsilon>0$ there exists $\delta>0$ such that whenever $\Delta=\{(s_i, [t_{i-1}, t_i]):1\leq i\leq N\}$ is a tagged partition of $[a,b]$ with $|\Delta|<\delta$, one has
\begin{equation}\label{2.1}
  {\bf d}(f(\Delta), \bar x)<\varepsilon,
\end{equation}
where
\begin{equation}\label{2.2}
  f(\Delta):=\sum_{i=1}^N(t_i-t_{i-1})f(s_i).
\end{equation}
In this case, we say that $\bar x$ is the Riemann integral of $f$ on $[a,b]$ and it is denoted as $\bar x=\int_a^bf(t)dt$.
\end{definition}

\begin{definition}
Let $X$ be a vector space equipped with a complete and invariant metric $\bf d$ and $\varphi:[a,b]\rightarrow X$. The function $\varphi$ is said to be differentiable at $t\in [a,b]$ if there exists a vector $x_t\in X$ such that
\begin{equation}\label{2-3}
  {\bf d}(0, \varphi(t+s)-\varphi(t)- sx_t)=o(s).
\end{equation}
We denote it by $\varphi'(t)=x_t$.

\end{definition}

It is not hard to verify that a function Riemann integrable must be bounded. We first obtain a characterization for Riemann integration via the following proposition. We give its proof for the sake of completeness.
\begin{pro}
Let $X$ be a vector space equipped with a complete and invariant metric $\bf d$ such that
\begin{equation}\label{2.3}
  {\bf d}(\lambda x,\lambda y)\leq \lambda{\bf d}(x,y)\ \ \lambda\in [0,1) \ {\it and} \ \forall x,y\in X.
\end{equation}
Let $f:[a,b]\rightarrow X$ be a function. Then $f$ is Riemann integration if and only if there exists a vector $\bar x\in X$ with the following property: for any $\varepsilon>0$ there exists a partition $\Delta_{\varepsilon}$ such that ${\bf d}(f(\Delta), \bar x)<\varepsilon$ whenever $\Delta=\{(s_i, [t_{i-1}, t_i]):1\leq i\leq N\}$ is a tagged partition of $[a,b]$ that refines $\Delta_{\varepsilon}$.
\end{pro}

{\it Proof.} The necessity part follows from the Riemann integration in Definition 2.2.

The sufficiency part. Let $\bar x\in X$ satisfy the corresponding property. Then one can check $f$ is bounded and let $L>0$ be a bound for the range of $f$ on $[a,b]$. Let $\varepsilon>0$ and take a partition $\Delta_{\varepsilon}=\{t_i:0\leq i\leq N\}$ of $[a,b]$ such that
\begin{equation}\label{2.4}
  {\bf d}(f(\Delta), \bar x)<\frac{\varepsilon}{2}
\end{equation}
whenever $\Delta$ is a tagged partition of $[a,b]$ that refines $\Delta_{\varepsilon}$.

Choose $\delta:=\min\{\frac{\varepsilon}{4LN}, 1\}$ and take any tagged partition $\Delta$ of $[a,b]$ with $|\Delta|<\delta$. We next prove that
\begin{equation}\label{2.5}
  {\bf d}(f(\Delta),\bar x)<\varepsilon.
\end{equation}
Put points of $\Delta$ and $\Delta_{\varepsilon}$ together and yield a new partition $\Delta_1$. We take the tags of $\Delta_1$ as follows:
\begin{itemize}
\item[1$^\circ$] the tag of each interval of $\Delta_1$ that coincides with an interval of $\Delta$ is the same as the tag of $\Delta$;
\item[2$^\circ$] the tags of $\Delta_1$ for the remaining intervals can be taken arbitrarily.
\end{itemize}
Without loss of generality, we can assume that $\{[c_k,d_k]:1\leq k\leq K\}$ are the intervals of $\Delta$ that contain points of $\Delta_{\varepsilon}$, saying $\{u_i^k:1\leq i\leq n_k-1\}\subset (c_k,d_k)$. Then $K\leq N$ and
\begin{equation}\label{2.6}
  c_k=u_0^k<u_1^k<\cdots<u^k_{n_k-1}<u^k_{n_k}=d_k.
\end{equation}
Suppose that $s_k$ is the tag of $\Delta$ for $[c_k,d_k]$ and $v_i^k$ is the tag of $\Delta_1$ for $[u_{i-1}^k, u_i^k]$ for each $i$. Then by \eqref{2.6}, one has
\begin{eqnarray*}
f(\Delta)-f(\Delta_1)&=&\sum_{k=1}^K\left((d_k-c_k)f(s_k)-\sum_{i=1}^{n_k}(u_i^k-u_{i-1}^k)f(v_i^k)\right)\\
&=&\sum_{k=1}^K\sum_{i=1}^{n_k}(u_i^k-u_{i-1}^k)(f(s_k)-f(v_i^k)).
\end{eqnarray*}
This and \eqref{2.3} imply that
\begin{eqnarray*}
{\bf d}(f(\Delta), f(\Delta_1))&\leq&\sum_{k=1}^K\sum_{i=1}^{n_k}{\bf d}((u_i^k-u_{i-1}^k)f(s_k), (u_i^k-u_{i-1}^k)f(v_i^k))\\
&\leq&\sum_{k=1}^K\sum_{i=1}^{n_k}(u_i^k-u_{i-1}^k){\bf d}(f(s_k), f(v_i^k))\\
&\leq&2L\sum_{k=1}^K(d_k-c_k)\\
&\leq& 2LN\delta\\
&<&\frac{\varepsilon}{2}
\end{eqnarray*}
as $L$ is the bound for the range of $f$ on $[a,b]$. Since $\Delta_1$ refines $\Delta_{\varepsilon}$, it follows from \eqref{2.4} that
\begin{eqnarray*}
{\bf d}(f(\Delta), \bar x)\leq {\bf d}(f(\Delta), f(\Delta_1))+ {\bf d}(f(\Delta_1), \bar x)<\frac{\varepsilon}{2}+\frac{\varepsilon}{2}=\varepsilon.
\end{eqnarray*}
This means that \eqref{2.5} holds. The proof is complete.\hfill$\Box$\\

The following proposition presents several equivalent criteria for Riemann integration whose proof follows from the Riemann integration in Definition 2.2 and Proposition 2.1.
\begin{pro}
Let $X$ be a vector space equipped with a complete and invariant metric $\bf d$ such that \eqref{2.3} holds and $f:[a,b]\rightarrow X$ be a function. Then the following statements are equivalent:
\begin{itemize}
\item[\rm(i)] The function $f$ is Riemann integrable on $[a,b]$;
\item[\rm(ii)] for any $\varepsilon>0$ there exists $\delta>0$ such that ${\bf d}(f(\Delta_1), f(\Delta_2))<\varepsilon$ for all tagged partition $\Delta_1$ and $\Delta_2$ of $[a,b]$ satisfying $\max\{|\Delta_1|, |\Delta_2|\}<\delta$;
\item[\rm(iii)] for any $\varepsilon>0$ there exists a partition $\Delta_{\varepsilon}$ such that ${\bf d}(f(\Delta_1),f(\Delta_2))<\varepsilon$ for all tagged partition $\Delta_1$ and $\Delta_2$ of $[a,b]$ that refine $\Delta_{\varepsilon}$;
\item[\rm(iv)] for any $\varepsilon>0$ there exists a partition $\Delta_{\varepsilon}$ such that ${\bf d}(f(\Delta_1), f(\Delta_2))<\varepsilon$ for all tagged partition $\Delta_1$ and $\Delta_2$ of $[a,b]$ that have the same points as $\Delta_{\varepsilon}$.
\end{itemize}
\end{pro}

{\it Proof.} (iv)$\Rightarrow$(iii): Let $\varepsilon>0$. Then there exists a partition $\Delta_{\varepsilon}=\{t_i:1\leq i\leq N\}$ such that ${\bf d}(f(\Delta_1),f(\Delta_2))<\frac{\varepsilon}{2}$ for all tagged partition $\Delta_1$ and $\Delta_2$ of $[a,b]$ that have the same points as $\Delta_{\varepsilon}$. Let $\Delta_0:=\{(t_i,[t_{i-1},t_i]):1\leq i\leq N\}$ be the tagged partition of $[a,b]$. For each $i$, let $A_i:=\{(t_i-t_{i-1})f(t): t\in [t_{i-1},t_i]\}$ and $A:=\sum_{i=1}^NA_i$. We claim that
\begin{equation}\label{2.7}
  {\bf d}(0, x)<\frac{\varepsilon}{2},\ \ \forall x\in {\rm co}(A-A).
\end{equation}
Indeed, let $x\in {\rm co}(A-A)$. Then there exist $\lambda_j\in [0,1]$ and $y_j,z_j\in A,j=1,\cdots,m$ such that
$$
\sum_{j=1}^m\lambda_j=1\ \ {\rm and} \ \ x=\sum_{j=1}^m\lambda_j(y_j-z_j).
$$
By virtue of \eqref{2.3}, one can yield that
\begin{eqnarray*}
{\bf d}(0,x)&=&{\bf d}(0,\sum_{j=1}^m\lambda_j(y_j-z_j))={\bf d}(\sum_{j=1}^m\lambda_jy_j, \sum_{j=1}^m\lambda_jz_j))\\
&\leq&{\bf d}\Big(\sum_{j=1}^m\lambda_jy_j, \lambda_1z_1+\sum_{j=2}^m\lambda_jy_j)\Big)+{\bf d}\Big(\lambda_1z_1+\sum_{j=2}^m\lambda_jy_j, \sum_{j=1}^m\lambda_jz_j\Big)\\
&=&{\bf d}(\lambda_1y_1, \lambda_1z_1)+{\bf d}(\sum_{j=2}^m\lambda_jy_j, \sum_{j=2}^m\lambda_jz_j)\\
&\leq&\lambda_1{\bf d}(y_1, z_1)+{\bf d}\Big(\sum_{j=2}^m\lambda_jy_j, \sum_{j=2}^m\lambda_jz_j\Big)\\
&\leq&\sum_{j=1}^m\lambda_j{\bf d}(y_j, z_j)<\sum_{j=1}^m\lambda_j\frac{\varepsilon}{2}=\frac{\varepsilon}{2}.
\end{eqnarray*}

Let $\Delta:=\{(v_k,[u_{k-1}, u_k]):1\leq k\leq K\}$ be a tagged partition of $[a,b]$ that refines $\Delta_{\varepsilon}$. Then for each $i$, we can assume that $u_{n_i}=t_i$ and consequently
$$
t_{i-1}=u_{n_{i-1}}<u_{n_{i-1}+1}<\cdots<u_{n_i-1}<u_{n_i}=t_i
$$
and
\begin{eqnarray*}
f(\Delta_0)-f(\Delta)&=&\sum_{i=1}^N\Big((t_i-t_{i-1})f(t_i)-\sum_{k=n_{i-1}+1}^{n_i}(u_k-u_{k-1})f(v_k)\Big)\\
&=&\sum_{i=1}^N\sum_{k=n_{i-1}+1}^{n_i}\frac{u_k-u_{k-1}}{t_i-t_{i-1}}((t_i-t_{i-1})f(t_i)-(t_i-t_{i-1})f(v_k))\\
&\in&\sum_{i=1}^N{\rm co}(A_i-A_i)\subset {\rm co}(A-A).
\end{eqnarray*}
This and \eqref{2.7} imply that
\begin{eqnarray*}
{\bf d}(f(\Delta_0),f(\Delta))={\bf d}(f(\Delta_0)-f(\Delta),0)<\frac{\varepsilon}{2}.
\end{eqnarray*}
Hence for any tagged partition $\Delta_1$ and $\Delta_2$ of $[a,b]$ that refine $\Delta_{\varepsilon}$, one has
$$
{\bf d}(f(\Delta_1),f(\Delta_2))\leq {\bf d}(f(\Delta_1),f(\Delta_0))+{\bf d}(f(\Delta_0),f(\Delta_2))<\frac{\varepsilon}{2}+\frac{\varepsilon}{2}=\varepsilon.
$$
This means that (iii) holds. The proof is complete.\hfill$\Box$\\

\noindent{\bf Remark 2.1.} For the case that $X$ is a Banach space, the metric can be given by the norm and this metric is complete and invariant and \eqref{2.3} is satisfied. Then Propositions 2.1 and 2.2 are valid and reduce to \cite[Theorem 3]{Go} and \cite[Theorem 5]{Go}, respectively. \hfill$\Box$\\

The following proposition provides a sufficient condition for Riemann integration in metrizable vector spaces.

\begin{pro}
Let $X$ be a vector space equipped with a complete and invariant metric $\bf d$ such that \eqref{2.3} holds and $f:[a,b]\rightarrow X$ be a function. Suppose that
\begin{equation}\label{2.8}
  \sup{\bf d}\Big(0,\sum_i(f(d_i)-f(c_i))\Big)
\end{equation}
is finite where the supremum is taken over all finite collections $\{[c_i,d_i]\}$ of nonoverlapping intervals in $[a,b]$. Then $f$ is Riemann integrable on $[a,b]$.
\end{pro}

{\it Proof.} Let $\varepsilon\in (0, 1)$. Suppose that $L>0$ is the bound of supremum in \eqref{2.8} taken over all finite collections of nonoverlapping intervals in $[a,b]$. Take $N\in\mathbb{N}$ such that $\frac{L}{N}(b-a)<\varepsilon$ and let $\Delta_{\varepsilon}:=\{t_i=a+\frac{i}{N}(b-a):0\leq i\leq N\}$ be a partition of $[a,b]$. Suppose that $\Delta_{1}:=\{(u_i, [t_{i-1}, t_i]): 1\leq i\leq N\}$ and $\Delta_{2}:=\{(v_i, [t_{i-1}, t_i]): 1\leq i\leq N\}$ are two tagged partition of $[a,b]$. Then partitions $\Delta_{1}$ and $\Delta_{2}$ have the same points as $\Delta_{\varepsilon}$ and
\begin{eqnarray*}
f(\Delta_1)-f(\Delta_2)=\sum_{i=1}^N(t_i-t_{i-1})(f(u_i)-f(v_i))=\frac{b-a}{N}\sum_{i=1}^N(f(u_i)-f(v_i)).
\end{eqnarray*}
Thus
\begin{eqnarray*}
{\bf d}(f(\Delta_1),f(\Delta_2))&=&{\bf d}(0, f(\Delta_1)-f(\Delta_2))\\
&\leq&\frac{b-a}{N}{\bf d}\Big(0, \sum_{i=1}^N(f(u_i)-f(v_i))\Big)\\
&\leq&\frac{L}{N}(b-a)<\varepsilon.
\end{eqnarray*}
Then Proposition 2.1(iv) implies that $f$ is Riemann integration on $[a,b]$. The proof is complete.\hfill$\Box$\\

The following theorem shows the existence of primitives on Riemann integration in metrizable vector spaces.
\begin{them}
Let $X$ be a vector space equipped with a complete and invariant metric $\bf d$ such that \eqref{2.3} holds and $f:[a,b]\rightarrow X$ be Riemann integrable on $[a,b]$. Define $F(t):=\int_a^tf(s)ds$ for any $t\in [a,b]$. Then for any $t$ at which $f$ is continuous, $F$ is differentiable and $F'(t)=f(t)$.
\end{them}

{\it Proof.} Let $t\in [a,b]$ be such that $f$ is continuous at $t$ and $\varepsilon>0$. Then there exists $r>0$ such that
\begin{equation}\label{2.10}
  {\bf d}(f(t+s), f(t))<\frac{\varepsilon}{2},\ \ \forall s\in (-r, r).
\end{equation}
Suppose that $s\in\mathbb{R}$ with $|s|<r$. Without loss of generality, we can assume that $s>0$. Then there exists $\delta>0$ such that for any tagged partitions $\Delta_1$ of $[a,t+s]$ and $\Delta_2$ of $[a,t]$ with $\max\{|\Delta_1|, |\Delta_2|\}<\delta$, one has
\begin{equation}\label{2.11}
  {\bf d}(F(t+s), f(\Delta_1))<\frac{s}{4}\varepsilon\ \ {\rm and} \ \ {\bf d}(F(t), f(\Delta_2))<\frac{s}{4}\varepsilon.
\end{equation}
Choose $N,K\in\mathbb{N}$ such that $\frac{t-a}{N}<\delta$ and $\frac{s}{K}<\delta$. Let $\Delta_1:=\{(t_i,[t_{i-1}, t_i]):1\leq i\leq N\}$ be a tagged partition of $[a,t]$ with $t_i:=a+\frac{i}{N}(t-a)(i=0,1,\cdots,N)$ and $\Delta_2:=\{(s_k,[s_{k-1}, s_k]):1\leq k\leq K\}$ be a tagged partition of $[t,t+s]$ with $s_k:=t+\frac{k}{K}s(k=0,1,\cdots,K)$ . Then $\Delta_{\varepsilon}:=\Delta_1\cup \Delta_2$ is a tagged partition of $[a,t+s]$ and $|\Delta_{\varepsilon}|<\delta$. By virtue of \eqref{2.3}, \eqref{2.10} and \eqref{2.11}, one has
\begin{eqnarray*}
&&{\bf d}(0, F(t+s)-F(t)-sf(t))\\
&=&{\bf d}\Big(0, F(t+s)-f(\Delta_{\varepsilon})+f(\Delta_{1})-F(t)+\sum_{k=1}^K(s_k-s_{k-1})(f(s_k)-f(t))\Big)\\
&\leq&{\bf d}(0, F(t+s)-f(\Delta_{\varepsilon}))+{\bf d}(0,f(\Delta_{1})-F(t))+{\bf d}\Big(0, \sum_{k=1}^K(s_k-s_{k-1})(f(s_k)-f(t))\Big)\\
&\leq&\frac{s}{4}\varepsilon+\frac{s}{4}\varepsilon+\sum_{k=1}^K\frac{s}{K}{\bf d}(0,f(s_k)-f(t))\\
&<&\frac{s}{2}\varepsilon+\sum_{k=1}^K\frac{s}{K}\frac{1}{2}\varepsilon=\varepsilon s.
\end{eqnarray*}
This means that
\begin{eqnarray*}
{\bf d}(0, F(t+s)-F(t)-sf(t))=o(s)
\end{eqnarray*}
and consequently $F'(t)=f(t)$. The proof is complete.\hfill$\Box$\\

The following result provides the theorem of calculus on Riemann integration in metrizable vector spaces.
\begin{them}
Let $X$ be a vector space equipped with a complete and invariant metric $\bf d$ such that \eqref{2.3} holds and $F:[a,b]\rightarrow X$ be differentiable on $[a,b]$. If $F'$ is continuous on $[a,b]$, then
\begin{equation}\label{2.12}
  \int_a^{\tau}F'(s)ds=F(\tau)-F(a),\ \ \forall \tau\in [a,b].
\end{equation}
\end{them}

{\it Proof.} Let $\tau\in [a,b]$ and $\varepsilon>0$. Noting that $F'$ is continuous on $[a,b]$, it follows that $\int_a^{\tau}F'(s)ds$ exists and it is denoted by $x_{\tau}$. Then there exists $\delta>0$ such that
\begin{equation}\label{2.13}
  {\bf d}(0,x_{\tau}-F'(\Delta))<\frac{\varepsilon}{2}
\end{equation}
holds for any tagged partition $\Delta$ of $[a,\tau]$ with $|\Delta|<\delta$.

Let $t\in [a,\tau]$. By the definition of $F'(t)$, there exists $r_t\in(0,\delta)$ such that
\begin{equation}\label{2.14}
  {\bf d}(0, F(t+s)-F(t)-sF'(t))<\frac{|s|}{2(b-a)}\varepsilon,\ \ \forall s\in (-r_t,r_t).
\end{equation}
Since $[a,b]$ is compact, there exist $t_1,\cdots,t_m\in [a,b]$ such that
\begin{equation}\label{2.15}
  [a,b]\subset\bigcup_{i=1}^m(t_i-\frac{1}{2}r_{t_i}, t_i+\frac{1}{2}r_{t_i}).
\end{equation}
Without loss of generality, we can assume that $a<t_1<\cdots<t_m<\tau$. By \eqref{2.15}, we can choose
\begin{equation}\label{2.16}
  u_i\in(t_i, t_i+\frac{1}{2}r_{t_i})\cap (t_i-\frac{1}{2}r_{t_i}, t_{i+1}), i=1,\cdots,m-1
\end{equation}
and  $u_0:=a,u_m:=\tau$. Let $\Delta_m:=\{(t_i, [u_{i-1}, u_i]):1\leq i\leq m\}$. Then $\Delta_m$ is the tagged partition of $[a,\tau]$ with $|\Delta_m|<\delta$ as $u_i-u_{i-1}=u_i-t_i+t_i-u_{i-1}<r_{t_i}<\delta$ for any $i$. Note that
\begin{eqnarray*}
&&(u_i-u_{i-1})F'(t_i)-F(u_i)+F(u_{i-1})\\
&=&(u_i-t_{i})F'(t_i)-F(u_i)+F(t_{i})+(t_i-u_{i-1})F'(t_i)-F(t_i)+F(u_{i-1}).
\end{eqnarray*}
By virtue of \eqref{2.14} and \eqref{2.16}, one has
\begin{eqnarray*}
 && {\bf d}(0,F'(t_i)(u_i-u_{i-1})-F(u_i)+F(u_{i-1}))\\
&\leq&{\bf d}(0, (u_i-t_{i})F'(t_i)-F(u_i)+F(t_{i}))+{\bf d}(0, (t_i-u_{i-1})F'(t_i)-F(t_i)+F(u_{i-1}))\\
&<&\frac{\varepsilon}{2(b-a)}(u_i-t_i)+\frac{\varepsilon}{2(b-a)}(t_i-u_{i-1})\\
&=&\frac{\varepsilon}{2(b-a)}(u_i-u_{i-1}).
\end{eqnarray*}
Note that
\begin{eqnarray*}
&&x_{\tau}-(F(\tau)-F(a))=x_{\tau}-F'(\Delta_m)+F'(\Delta_m)-(F(\tau)-F(a))\\
&=&x_{\tau}-F'(\Delta_m)+\sum_{i=1}^m((u_i-u_{i-1})F'(t_i)-F(u_i)+F(u_{i-1})).
\end{eqnarray*}
Therefore
\begin{eqnarray*}
 && {\bf d}(0, x_{\tau}-(F(\tau)-F(a)))\\
&\leq&{\bf d}(0, x_{\tau}-F'(\Delta_m))+\sum_{i=1}^m{\bf d}(0, (u_i-u_{i-1})F'(t_i)-F(u_i)+F(u_{i-1}))\\
&<&\frac{\varepsilon}{2}+\sum_{i=1}^m\frac{\varepsilon}{2(b-a)}(u_i-u_{i-1})\\
&<&\varepsilon.
\end{eqnarray*}
By taking the limit as $\varepsilon\rightarrow 0^+$, one can yield ${\bf d}(0, x_{\tau}-(F(\tau)-F(a)))=0$ and thus $x_{\tau}=F(\tau)-F(a)$. This means that \eqref{2.12} holds. The proof is complete.\hfill$\Box$

\setcounter{equation}{0}
\section{Main Results}

In this section, we study the relationship between continuity and Riemann integrability of a function in metrizable vector spaces. As main results, we prove that $l^1(\Gamma)$ ($\Gamma$ uncountable) and $\mathbb{R}^{\omega}$ equipped with product topology are metrizable vector spaces having the Lebesgue property (see Definition 3.2). For convenience to study Riemann integrability and continuity of functions, we introduce some notations.

Let $X$ be a vector space with a complete and invariant metric $\bf d$ and $f:[a,b]\rightarrow X$. For any closed subinterval $I$ in $[a,b]$, denote by
\begin{equation*}
  \omega(f,I):=\sup\{{\bf d}(f(u), f(v)): u,v\in I\}.
\end{equation*}
the oscillation of $f$ on the interval $I$. For each partition $\Delta=\{t_i:0\leq i\leq N\}$ of $[a,b]$, let
\begin{equation*}
  \omega(f,\Delta):=\sum_{i=1}^N\omega(f,[t_{i-1}, t_i])(t_i-t_{i-1}).
\end{equation*}
For each $t\in (a,b)$, denote by
\begin{equation*}
  \omega(f,t):=\lim_{\delta\rightarrow 0^+}\omega(f,[t-\delta, t+\delta]).
\end{equation*}
the oscillation of $f$ at $t$ and let
\begin{equation*}
  \omega(f,a):=\lim_{\delta\rightarrow 0^+}\omega(f,[a, a+\delta])\ {\rm and} \
  \omega(f,b):=\lim_{\delta\rightarrow 0^+}\omega(f,[b-\delta, b]).
\end{equation*}
It is easy to verify that $f$ is continuous at $t\in [a,b]$ if and only if $\omega(f,t)=0$, and the set $\{t\in [a,b]:\omega(f,t)\geq r\}$ is closed for any $r\in\mathbb{R}$.\\

\begin{definition}
Let $X$ be a vector space equipped with a complete and invariant metric ${\bf d}$ and $f:[a,b]\rightarrow X$ be a function. The function is said to be Darboux integrable if for any $\varepsilon>0$ there exists a partition $\Delta_{\varepsilon}$ of $[a,b]$ such that $\omega(f,\Delta)<\varepsilon$ whenever $\Delta$ is a partition of $[a,b]$ that refines $\Delta_{\varepsilon}$.
\end{definition}

The following proposition is on the continuity almost everywhere of functions in metric vector spaces. The proof is inspired by that of \cite[Theorem 18]{Go} and we present it in detail for the sake of completeness.
\begin{pro}
Let $X$ be a vector space equipped with a complete and invariant metric $\bf d$ and $f:[a,b]\rightarrow X$ be a function. Then $f$ is Darboux integrable if and only if $f$ is bounded and continuous almost everywhere on $[a,b]$.
\end{pro}

{\it Proof.} The necessity part. It is easy to verify that $f$ is bounded. For each $r>0$, let $E_r:=\{t\in [a,b]:\omega(f,t)\geq r\}$ and let
\begin{equation}\label{3.1a}
  E:=\{t\in [a,b]:\omega(f,t)>0\}.
\end{equation}
Then $E=\bigcup_{r>0}E_r$. To complete the proof of the necessity part, we only need to show that $m(E)=0$.

Suppose on the contrary that $m(E)>0$. Then there exists $r>0$ such that $m(E_r)>0$ as $E=\bigcup_{r>0}E_r$. Let $\Delta$ be any partition of $[a,b]$ and $\Delta_1$ be a partition of $[a,b]$ refining $\Delta$. Denote by $\{I_j:1\leq j\leq p\}$ the all closed intervals of $\Delta_1$ such that $I_j\cap E_r\not=\emptyset$. Then
$$
m(E_r)=m\big(\bigcup_{j=1}^p(I_j\cap E_r)\big)\leq\sum_{j=1}^p m(I_j\cap E_r)\leq\sum_{j=1}^p m(I_j).
$$
This implies that
$$
\omega(f,\Delta_1)\geq\sum_{j=1}^p\omega(f,I_j)m(I_j)\geq r\sum_{j=1}^p m(I_j)\geq rm(E_r)>0,
$$
which is a contradiction to the Darboux integrability of $f$.

The suffiency part. Let $L>0$ be a bound for $f$ and $\varepsilon>0$. Choose $N\in\mathbb{N}$ such that $\frac{b-a}{N}<\frac{\varepsilon}{2}$. Let
$$
E_N:=\left\{t\in [a,b]: \omega(f,t)\geq\frac{1}{N}\right\}.
$$
Then $m(E_N)=0$ as $f$ is continuous almost everywhere on $[a,b]$. Since $E_N\subset [a,b]$ is closed and bounded, it follows that there exist finite open intervals, saying $\{(c_i,d_i):1\leq i\leq p\}$, in $[a,b]$ such that
\begin{equation}\label{3.2a}
d_i<c_{i+1}, i=1,\cdots,p-1,\ \   E_N\subset\bigcup_{i=1}^p(c_i,d_i)\ \ {\rm and}\ \ \sum_{i=1}^p(d_i-c_i)<\frac{\varepsilon}{4L}.
\end{equation}
Let $\Delta_{\varepsilon}^1:=\{[c_i,d_i]:1\leq i\leq p\}$.

Let $[u,v]$ be an interval in $[a,b]$ contiguous to the intervals of $\Delta_{\varepsilon}^1$ such that $[u,v]\cap E_N=\emptyset$. Then for any $t\in [u,v]$, there exists $\delta_t>0$ such that $\omega(f, [t-\delta_t,t+\delta_t])<\frac{1}{N}$. This implies that the collection $\{([t-\delta_t,t+\delta_t):t\in [u,v]\}$ is an open cover of $[u,v]$ and thus there exists a finite subcover. Choose endpoints of the intervals from the finite subcover that belong to $(u,v)$ together with $u$ and $v$ form a partition of $[u,v]$. Do this process for all intervals in $[a,b]$ contiguous to the intervals of $\Delta_{\varepsilon}^1$ that  do not intersect $E_N$.  Denote by $\Delta_{\varepsilon}^2$ all of these intervals. Then the intervals of $\Delta_{\varepsilon}^1$ and $\Delta_{\varepsilon}^2$ form a partition of $[a,b]$ and it is denoted by $\Delta_{\varepsilon}$.

Now let $\Delta$ be a partition of $[a,b]$ that refines $\Delta_{\varepsilon}$. Let $\Delta^1$ and $\Delta^2$ be the intervals of $\Delta$ that are entirely contained within intervals of $\Delta_{\varepsilon}^1$ and $\Delta_{\varepsilon}^2$, respectively. Then by \eqref{3.2a}, one has
\begin{eqnarray*}
\omega(f,\Delta)&=&\sum_{I\in\Delta^1}\omega(f,I)m(I)+\sum_{I\in\Delta^2}\omega(f,I)m(I)\\
&\leq&2L\cdot\frac{\varepsilon}{4L}+\frac{1}{N}(b-a)\\
&<&\frac{\varepsilon}{2}+\frac{\varepsilon}{2}=\varepsilon.
\end{eqnarray*}
This means that $f$ is Darboux integrable on $[a,b]$. The proof is complete.\hfill$\Box$\\

The following corollary follows from Propositions 2.2 and 3.1.
\begin{coro}
Let $X$ be a vector space equipped with a complete and invariant metric $\bf d$ such that \eqref{2.3} holds and $f:[a,b]\rightarrow X$ be a function. If $f$ is continuous almost everywhere on $[a,b]$, then $f$ is Riemann integrable on $[a,b]$.
\end{coro}

From Corollary 3.1, it is an interesting problem to determine metrizable vector spaces in which every Riemann integrable function is continuous almost everywhere. This is equivalent to determine spaces with the property that Riemann integrability of functions is equivalent to Darboux integrability. This property is known as the Lebesgue property  since Lebesgue first proved that $\mathbb{R}$ has such property. In the context of metrizable vector space, we consider the following definition.

\begin{definition}
Let $X$ be a vector space with a complete and invariant metric $\bf d$. We say that $X$ has the Lebesgue property  if for every Riemann integrable function $f:[a,b]\rightarrow X$ is continuous almost everywhere on $[a,b]$.
\end{definition}

The following proposition is a useful tool in determining whether a space has the Lebesgue property. The proof is easy to obtain and thus it is omitted.
\begin{pro}
Let $X$ be a vector space with a complete and invariant metric $\bf d$ and $Y$ be a subspace of $X$.
\begin{itemize}
\item[\rm(i)] If $X$ has the Lebesgue property, then $Y$ has the Lebesgue property.
\item[\rm(ii)] If $Y$ does not have the Lebesgue property, then $X$ does not have the Lebesgue property.
\end{itemize}
\end{pro}

In the context of Banach spaces, Lebesgue first proved that $\mathbb{R}$ has the property in Definition 3.2. By arguing in each coordinate, every finite Euclidean space has the Lebesgue property. The following two examples show that $l^p(1\leq p\leq+\infty)$ do not have the Lebesgue property. \\

\noindent{\bf Example 3.1.} Let $p\in(1,+\infty)$ and $\{r_n\}$ be all rational numbers in $[0,1]$. Define $f:[0,1]\rightarrow l^p$ as $f(t)=0$ if $t$ is irrational and $f(t)=e_n$ if $t=r_n$. Then $f$ is Riemann integrable on $[0,1]$ since it is of outside bounded on $[0,1]$ (see \cite[Definition 6 and Theorem 9]{Go}). However, it is clear that $f$ is not continuous at all irrational numbers and thus not continuous almost everywhere on $[0,1]$. Furthermore, this conclusion on $f$ shows that $c_0$ does not have the Lebesgue property .\hfill$\Box$\\


\noindent{\bf Example 3.2.} Let $f:[0,1]\rightarrow l^{\infty}$ be defined as $f(t):=(c_k(t))_{k=1}^{\infty}\in l^{\infty}$ for any
$$t=\sum_{k=1}^{\infty}\frac{c_k(t)}{2^k}\in [0,1]$$
where each of the digits $c_k(t)$ is $0$ or $1$.  Then one can verify that $f$ is Riemann integrable on $[0,1]$ by Proposition 2.2. However, $f$ is not
continuous at any $t\in [0,1]$. Hence $l^{\infty}$ does not have the Lebesgue property . \hfill$\Box$\\

Further, da Rocha \cite{Ro} proved that infinite dimensional and uniformly convex Banach spaces do not have the Lebesgue property  and consequently neither do infinite dimensional Hilbert spaces and $L^p[a,b](1\leq p\leq\infty)$. Therefore, to seek spaces with the Lebesgue property  is a natural and interesting problem. It is known that Nemirovski, Ochan and Rejouani \cite{NOR} and da Rocha \cite{Ro} proved that the space $l^1$ has the Lebesgue property, and da Rocha \cite{Ro} further proved that the Tsirelson space has the Lebesgue property. In this section, we further prove that $l^1(\Gamma)$($\Gamma$ uncountable) also possesses the Lebesgue property. We first recall the definition of $l^1(\Gamma)$($\Gamma$ uncountable).

Let $\Gamma$ be an uncountable index set. The space $l^1(\Gamma)$ consists of all elements that vanish at all but a countable number of members of $\Gamma$ and are absolutely summable; that is,
$$
x=(x_{\alpha})_{\alpha\in\Gamma}\Leftrightarrow \{\alpha\in\Gamma: x_{\alpha}\not=0\}\ {\it is\ countable\ and} \ \sum_{\alpha\in\Gamma}|x_{\alpha}|<+\infty.
$$
Let $\|x\|:=\sum_{\alpha\in\Gamma}|x_{\alpha}|$ for any $x=(x_{\alpha})_{\alpha\in\Gamma}\in l^1(\Gamma)$. It is easy to verify that $l^1(\Gamma)$ is a Banach space and consequently a metrizable vector space.\\

\noindent{\bf Remark 3.1.} The Banach space $l^1(\Gamma)$($\Gamma$ uncountable) is different from $l^1$, since $l^1$ is separable but $l^1(\Gamma)$ is not; and from \cite[Example 1.4]{Ph}, the norm in $l^1$ is G\^{a}teaux differentiable at those points $x=(x_n)$ for which $x_n\not=0$ for all $n$; while the norm in $l^1(\Gamma)$ is not G\^{a}teaux differentiable at any point.\hfill$\Box$

\begin{them}
The space $l^1(\Gamma)$ has the Lebesgue property .
\end{them}

{\it Proof.} Suppose that $f:[a,b]\rightarrow l^1(\Gamma)$ be a Riemann integrable function. Then the range of $f$ is a bounded subset of $l^1(\Gamma)$. For each $\alpha\in\Gamma$, denote by $\pi_{\alpha}(x):=x_{\alpha}$ the projection mapping on $l^1(\Gamma)$.

Suppose on the contrary that $f$ is not continuous almost everywhere on $[a,b]$. We next show that $f$ is not Riemann integrable, which is a  contradiction.

Let $$E:=\{t\in [a,b]: \omega(f,t)>0\}.$$
Then $E$ has positive measure and thus there exists $r>0$ such that $m(E_r)>0$ where
\begin{equation}\label{3-3}
E_r:=\{t\in [a,b]: \omega(f,t)\geq r\}.
\end{equation}
(Otherwise, $m(E_r)=0$ for all $r>0$ and consequently $E$ has measure zero as $E=\bigcup_{n=1}^{\infty}E_{\frac{1}{n}}$, which is a contradiction).

Let $\varepsilon_0:=\frac{rm(E_r)}{8}>0$ and $\delta>0$ be arbitrarily given. Choose $N\in\mathbb{N}$ such that $\frac{1}{N}<\delta$ and let
$$
\Delta_N:=\left\{\frac{i}{N}:0\leq i\leq N\right\}
$$
be a partition of $[0,1]$. Denote by $\{[c_k,d_k]:1\leq k\leq p\}$ the all intervals in $\Delta_N$ satisfying $m([c_k,d_k]\cap E_r)>0$. Then
\begin{equation}\label{3-4}
  m(E_r)=m\Big(\bigcup_{k=1}^p[c_k,d_k]\cap E_r\Big)\leq\sum_{k=1}^pm([c_k,d_k])=\frac{p}{N}.
\end{equation}
For each $\alpha\in\Gamma$, let $G_{\alpha}$ be the set of all discontinuous points of $\pi_{\alpha}\circ f$ on $[0,1]$. If $m(G_{\alpha})>0$ for some $\alpha\in\Gamma$, then $\pi_{\alpha}\circ f$ is not Riemann integrable as $\mathbb{R}$ has the Lebesgue property  and consequently $f$ is not Riemann integrable. We next consider the case that $m(G_{\alpha})=0$ for all $\alpha\in\Gamma$.

Let $\Gamma_0\subset \Gamma$ be a countable subset and $I_0=\emptyset$. Then $G_0:=\bigcup_{\alpha\in\Gamma_0}G_{\alpha}$ has measure zero. Note that $m([c_1,d_1]\cap E_r)>0$ and then we can take $u_1\in(E_r\backslash G_0)\cap (c_1,d_1)$ such that $\omega(f,u_1)\geq r$. Thus there exists $v_1\in (c_1,d_1)$ such that
$$
\|f(u_1)-f(v_1)\|\geq\frac{r}{2}.
$$
Let $x^1:=f(u_1)-f(v_1)=(x^1_{\alpha})_{\alpha\in\Gamma}$. Then there exists $\Gamma_1\subset\Gamma$ countable such that
\begin{equation}\label{3-5}
  x^1_{\alpha}=0, \forall \alpha\in\Gamma\backslash\Gamma_1.
\end{equation}
Take a finite index subset $I_1\subset\Gamma_1$ such that
\begin{equation}\label{3-6}
  \sum_{\alpha\in\Gamma_1\backslash I_1}|x^1_{\alpha}|<\frac{\varepsilon_0}{2}.
\end{equation}
Let $G_1:=\bigcup_{\alpha\in\Gamma_0\cup\Gamma_1}G_{\alpha}$. Then $m(G_1)=0$ and thus we can choose $u_2\in(E_r\backslash G_1)\cap(c_2,d_2)$ such that $\omega(f,u_2)\geq r$. Since $\pi_{\alpha}\circ f$ is continuous at $u_2$ for each $\alpha\in\Gamma_0\cup\Gamma_1$, then there exists $v_2\in(c_2,d_2)$ such that
\begin{equation}\label{3-7}
  \|f(u_2)-f(v_2)\|\geq\frac{r}{2}\ \ {\rm and}\ \ \sum_{\alpha\in I_1}|\pi_{\alpha}\circ f(v_2)-\pi_{\alpha}\circ f(u_2)|<\frac{\varepsilon_0}{2^{1}}.
\end{equation}
Let $x^2:=f(u_2)-f(v_2)=(x^2_{\alpha})_{\alpha\in\Gamma}$. Then there exists $\Gamma_2\subset\Gamma$ countable such that
\begin{equation}\label{3-8}
  x^2_{\alpha}=0, \forall \alpha\in\Gamma\backslash\Gamma_2.
\end{equation}
Take a finite index subset $I_2\subset\Gamma_2$ such that
\begin{equation*}
  \sum_{\alpha\in\Gamma_2\backslash I_2}|x^2_{\alpha}|<\frac{\varepsilon_0}{2^{2}}.
\end{equation*}
Let $G_2:=\bigcup\limits_{i=0}^2\bigcup\limits_{\alpha\in\Gamma_i}G_{\alpha}$. Then $m(G_2)=0$ and thus we can choose $u_3\in(E_r\backslash G_2)\cap(c_3,d_3)$ such that $\omega(f,u_3)\geq r$. Since $\pi_{\alpha}\circ f$ is continuous at $u_3$ for each $\alpha\in\bigcup\limits_{i=0}^2\Gamma_i$, then there exists $v_3\in(c_3,d_3)$ such that
\begin{equation*}
  \|f(u_3)-f(v_3)\|\geq\frac{r}{2}\ \ {\rm and}\ \ \sum_{\alpha\in I_1\cup I_2}|\pi_{\alpha}\circ f(v_3)-\pi_{\alpha}\circ f(u_3)|<\frac{\varepsilon_0}{2^{2}}.
\end{equation*}
Let $x^3:=f(u_3)-f(v_3)=(x^3_{\alpha})_{\alpha\in\Gamma}$. We continue this process for $p$ steps and get sequences $\{G_k,\Gamma_k,I_k,u_{k},v_{k},x^k):1\leq k\leq p\}$ in $[0,1]\times\Gamma\times\Gamma\times(E_r\backslash G_k)\cap(c_{k},d_{k})\times (c_{k},d_{k})\times l^1(\Gamma)$ such that
\begin{equation}\label{3-9}
\Gamma_k\subset \Gamma\ {\rm countable}, \ I_k\subset\Gamma_k\ {\rm finite}, \ G_k:=\bigcup\limits_{i=0}^k\bigcup\limits_{\alpha\in\Gamma_i}G_{\alpha}\  {\rm with} \ m(G_k)=0,
\end{equation}
\begin{equation}\label{3-10}
\omega(f,u_{k})\geq r, x^{k}=f(u_{k})-f(v_{k})=(x^{k}_{\alpha})_{\alpha\in\Gamma} \ \ {\rm with} \ \sum_{\alpha\in\Gamma_{k}\backslash I_k}|x^{k}_{\alpha}|<\frac{\varepsilon_0}{2^{k}}
\end{equation}
and
\begin{equation}\label{3-11}
\|f(u_{k})-f(v_{k})\|\geq\frac{r}{2}\ \ {\rm and}\ \ \sum_{i=1}^k\sum_{\alpha\in I_i}|\pi_{\alpha}\circ f(v_{k})-\pi_{\alpha}\circ f(u_{k})|<\frac{\varepsilon_0}{2^{k}}.
\end{equation}
For each $k:1\leq k\leq p$, let
$$
y^{k}:=\sum_{\alpha\in I_k\big\backslash \bigcup\limits_{i=1}^{k-1}I_{i}}x^k_{\alpha}e^{\alpha},
$$
where $e^{\alpha}$ is a vector in $l^1(\Gamma)$ satisfying $\pi_{\beta}(e^{\alpha})=1$ if $\beta=\alpha$ and $\pi_{\beta}(e^{\alpha})=0$ if $\beta\not=\alpha$.
By \eqref{3-10} and \eqref{3-11}, one has
\begin{equation}\label{3-12}
\|x^k-y^{k}\|=\sum_{i=1}^{k-1}\sum_{\alpha\in I_i}|x^k_{\alpha}|+\sum_{\alpha\in \Gamma_k\backslash I_{k}}|x^k_{\alpha}|<\frac{\varepsilon_0}{2^{k-1}}
\end{equation}
and
\begin{equation}\label{3-13}
  \|y^k\|=\|x^k\|-\|x^k-y^k\|\geq\frac{r}{2}-\frac{\varepsilon_0}{2^{k-1}}.
\end{equation}
This implies that
\begin{eqnarray*}
  \left\|\sum_{k=1}^px^k\right\|&\geq& \left\|\sum_{k=1}^py^k\right\|- \left\|\sum_{k=1}^p(y^k-x^k)\right\|\\
  &=&\sum_{k=1}^p\|y^k\|- \left\|\sum_{k=1}^p(y^k-x^k)\right\|\\
  &\geq&\sum_{k=1}^p\|y^k\|- \sum_{k=1}^p\|y^k-x^k\|\\
  &\geq&\sum_{k=1}^p\big(\frac{r}{2}-\frac{\varepsilon_0}{2^{k-1}}\big)-\sum_{k=1}^p\frac{\varepsilon_0}{2^{k-1}}\\
  &\geq&\frac{pr}{2}-2\varepsilon_0.
\end{eqnarray*}
Let $\Delta_1$ and $\Delta_2$ be two partitions of $[0,1]$ that have the same points as $\Delta_N$. The tags of $\Delta_1$ and $\Delta_2$ are $u_k$ and $v_k$ respectively in the intervals $[c_k,d_k]$ for $1\leq k\leq p$ and the tags of $\Delta_1$ and $\Delta_2$ are the same in the remaining intervals. By virtue of \eqref{3-4}, one has
\begin{eqnarray*}
\|f(\Delta_1)-f(\Delta_2)\|=\left\|\sum_{k=1}^p\frac{1}{N}x^k\right\|\geq\frac{pr}{2N}-\frac{2}{N}
\varepsilon_0\geq\frac{rm(E_r)}{2}-\frac{rm(E_r)}{4}=\frac{rm(E_r)}{4}>0.
\end{eqnarray*}
By virtue of Proposition 2.2, one can yield that $f$ is not Riemann integrable. The proof is complete.\hfill$\Box$\\

As one of main results in this paper, we discovery a metrizable vector space that has the Lebesgue property as said in Definition 3.2.

Let $\mathbb{R}^\omega$ be the countably infinite product of $\mathbb{R}$ with itself; that is,
\begin{equation}\label{3.1}
\mathbb{R}^\omega:=\prod_{i\in \mathbb{N}}X_i
\end{equation}
where $X_i=\mathbb{R}$ for each $i$. Let $\mathbb{R}^\omega$ be equipped with the product topology that is denoted by $\tau_{\pi}$. Then the product space $(\mathbb{R}^\omega,\mathcal{T}_{\pi})$ is a vector space. Further, let ${\bf d}:\mathbb{R}^\omega\times\mathbb{R}^\omega\rightarrow \mathbb{R}$ be defined as
\begin{equation}\label{3.2}
  {\bf d}(x,y):=\sup_{i\geq 1}\frac{\min\{|x_i-y_i|,1\}}{i}\ \ \forall x=(x_i)_{i\in\mathbb{N}}, y=(y_i)_{i\in\mathbb{N}}\in\mathbb{R}^\omega.
\end{equation}
It is easy to verify that $\bf d$ in \eqref{3.2} is a complete and variant metric and by \cite[Theorem 20.5]{M}, the metric $\bf d$ on $\mathbb{R}^\omega$ is compatible with the product topology $\mathcal{T}_{\pi}$. As one main result of this paper, we prove that $(\mathbb{R}^\omega,\bf d)$ has the Lebesgue property. Before this, we need the following lemma.

\begin{lem}
Let $\bf d$ be given as \eqref{3.2} and $f:[a,b]\rightarrow (\mathbb{R}^{\omega}, \bf d)$. Then $f$ is continuous if and only if $\pi_i\circ f:[a,b]\rightarrow \mathbb{R}$ is continuous for each $i$, where $\pi_j: \mathbb{R}^{\omega}\rightarrow X_j=\mathbb{R}$ is the projection mapping associate with the index $j$; that is $\pi_j((x_i)_{i\in\mathbb{N}})=x_j$.
\end{lem}

{\it Proof.} Note that the projection mapping $\pi_i$ is continuous and thus the necessity part holds.

The sufficiency part. Let $t\in [a,b]$ and $\varepsilon>0$. Choose $N\in\mathbb{N}$ such that
\begin{equation}\label{3-1}
  \sum_{i=N+1}^{\infty}\frac{1}{2^i}<\frac{\varepsilon}{2}.
\end{equation}
Since $\pi_i\circ f:[a,b]\rightarrow \mathbb{R}$ is continuous for each $i$, then there exists $\delta>0$ such that
\begin{equation}\label{3-2}
  \frac{|\pi_i(f(s))-\pi_i(f(t))|}{2^i}<\frac{\varepsilon}{2N}\ \ \forall s\in (t-\delta, t+\delta)\ {\rm and} \ \forall i=1,\cdots,N.
\end{equation}
Then for any $s\in (t-\delta, t+\delta)$, it follows from \eqref{3-1} and \eqref{3-2} that
\begin{eqnarray*}
{\bf d}(f(s),f(t))&=&\sum_{i=1}^{\infty} \frac{\min\{|\pi_i(f(s))-\pi_i(f(t))|,1\}}{2^i}\\
&\leq&\sum_{i=1}^N \frac{|\pi_i(f(s))-\pi_i(f(t))|}{2^i}+\sum_{i=N+1}^{\infty} \frac{1}{2^i}\\
&<&\sum_{i=1}^N \frac{\varepsilon}{2N}+\frac{\varepsilon}{2}=\varepsilon.
\end{eqnarray*}
Hence $f$ is continuous at $t$. The proof is complete.\hfill$\Box$

\begin{them}
The space $\mathbb{R}^{\omega}$ equipped with the product topology $\mathcal{T}_{\pi}$ has the Lebesgue property.
\end{them}

{\it Proof.} Let $f:[a,b]\rightarrow\mathbb{R}^{\omega}$ be a Riemann integrable function. Then one can verify that $\pi_i\circ f:[a,b]\rightarrow \mathbb{R}$ is Riemann integrable for each $i$. Since $\mathbb{R}$ has the Lebesgue property , then for each $i$, one has that $\pi_i\circ f$ is continuous almost everywhere on $[a,b]$ and
$$
E_i:=\{t\in [a,b]: \pi_i\circ f\ {\rm is\ discontinuous\ at} \ t\}
$$
has measure zero. Let $E:=\bigcup_{i=1}^{\infty}E_i$. Then the measure of $E$ is zero and it follows from Lemma 3.1 that the set with all continuous points $f$ in $[a,b]$ is $[a,b]\backslash E$. This implies that $f$ is continuous almost everywhere on $[a,b]$ and consequently $\mathbb{R}^{\omega}$ equipped with the product topology $\mathcal{T}_{\pi}$ has the Lebesgue property. The proof is complete.\hfill$\Box$\\

The following corollary follows from Proposition 3.2 and Theorem 3.2.
\begin{coro}
Let $\bf d$ be defined as in \eqref{3.2}. Then $l^p(1\leq p\leq+\infty)$, as subspaces of $(\mathbb{R}^{\omega},\bf d)$, have the Lebesgue property.
\end{coro}

\noindent{\bf Remark 3.2.} It is known from Examples 3.1 and 3.2 that $l^p(1\leq p\leq+\infty)$ are Banach spaces having no Lebesgue property. However, as said in Corollary 3.2, these spaces equipped with subspace topology of $(\mathbb{R}^{\omega},\bf d)$ possess the Lebesgue property.

\section{Conclusions}
This paper is devoted  to the study of Riemann integration and the Lebesgue property in the framework of metrizable vector spaces. Many of the real-valued results concerning the Riemann integration are still valid and we discovery that $\mathbb{R}^{\omega}$ equipped with product topology has the Lebesgue property. This enable us to consider the issue of Riemann integration and the Lebesgue property in the general topological vector spaces. Such issue may be more difficult and challenging. For example, let $\mathbb{R}^{\omega}$ in \eqref{3.1} be equipped with box topology $\mathcal{T}_B$ (different from product topology $\mathcal{T}_{\pi}$) and one can consider the Riemann integration of functions as Definition 2.2 via box topology. Consider a function $f:[a,b]\rightarrow (\mathbb{R}^{\omega}, \mathcal{T}_B)$ defined by $f(t):=(t,\cdots,t,\cdots,)$ for any $t\in [a,b]$. It is not hard to verify that $f$ is Riemann integrable with respect to the box topology whereas $f$ is not continuous on $[a,b]$. This means that the Lebesgue property of $\mathbb{R}^{\omega}$ may not remain valid in the sense of box topology.\\\\

\end{document}